\documentclass[12pt,reqno]{article}
\usepackage{amssymb,amscd,amsmath,amsthm,color}
\newtheorem{theorem}{Theorem}[section]

\newtheorem{cor}[theorem]{Corollary}

\usepackage{enumerate,amssymb}
\theoremstyle{definition}

\theoremstyle{remark}

\numberwithin{equation}{section}

\def\bP{\mathbb{P}}

\def\bM{\mathbb{M}}

\def\bR{\mathbb{R}}

\def\bR{\mathbb{R}}

\begin{document}
\baselineskip=15pt

\title{Convex maps on $\bR^n$ and positive definite matrices }

\author{  Jean-Christophe Bourin and  Jingjing Shao}

\date{ }

\maketitle

\vskip 10pt\noindent
{\small
{\bf Abstract.} We  obtain several convexity statements involving positive definite matrices. In  particular, if $A,B,X,Y$ are invertible matrices and $A,B$ are positive, we show that the map
$$
(s,t) \mapsto {\mathrm{Tr\,}} \log \left(X^*A^sX + Y^*B^tY\right)
$$
is jointly convex on $\bR^2$.  This is related to some exotic matrix H\"older inequalities such as
$$
  \left\| \sinh\left( \sum_{i=1}^m A_iB_i\right) \right\| \le  \left\| \sinh\left( \sum_{i=1}^m A_i^p\right) \right\|^{1/p} \left\| \sinh\left( \sum_{i=1}^m B_i^q\right) \right\|^{1/q}
$$
for all  positives matrices $A_i,  B_i$, such that $A_iB_i=B_iA_i$, conjugate exponents $p,q$ and unitarily invariant norms $\|\cdot\|$.
Our approach to obtain these results consists in studying the behaviour of some functionals along the geodesics of the Riemanian manifold of positive definite matrices.

}

\section{Convex and log-convex maps}

This short note aims to point out some convex maps involving positive definite matrices. 
We denote by $\bM_n$ the space of $n$-by-$n$ matrices with complex entries, and by $\bP_n$ its positive definite cone. A non-negative, continuous function $f(t)$ defined on $[0,\infty)$ is geometrically convex if $f(\sqrt{ab})\le \sqrt{f(a)f(b)}$ for all $a,b>0$, equivalently if $\log f(e^ t)$ is convex on $\bR$. Note that a  function $\varphi(t)$ on $(0,\infty)$ satifies the geometric-arithmetic convexity  inequality
$$
\varphi(\sqrt{ab}) \le \frac{\varphi(a)+\varphi(b)}{2}, \qquad a,b>0,
$$
if and only if $e^{ \varphi(t)}$ is geometrically convex, equivalently $\varphi(e^ t)$ is convex on $\bR$. This convexity property can be extended to the matrix setting as follows. 

\vskip 10pt
\begin{theorem}\label{th-trace} Let $\varphi(t)$ be a nondecreasing function defined on $(0,\infty)$ such that $ \varphi(e^t)$ is  convex.  Let $A_i\in\bP_n$ and $X_i\in\bM_n$ be invertible, $i=1,\ldots, m$. Then, the map
$$
(t_1,\ldots,t_m) \mapsto {\mathrm{Tr\,}}\varphi\left(\sum_{i=1}^m X_i^*A_i^{t_i} X_i\right)
$$
is jointly convex on $\bR^m$.
\end{theorem}

\vskip 10pt
Letting $\varphi(t)=\log t$, we get the statement of the Abstract. Theorem \ref{th-trace} can be derived from the following more general log-convexity theorem. Recall that a symmetric norm on $\bM_n$ satifies $\| UAV\|=\| A\|$ for all $A\in\bM_n$ and all unitary matrices $U,V\in\bM_n$. We denote by $\bM_n^+$ the positive semi-definite cone of $\bM_n$. A positive linear map $\Phi:\bM_n\mapsto \bM_d$ satifies $\Phi(\bM_n^+)\subset \bM_d^+$. A classical example  is the Schur multipler $A\mapsto Z\circ A$ with $Z\in\bM_n^+$.

\vskip 10pt
\begin{theorem}\label{th-norm}   Let $A_i\in\bM_n^+$ and $X_i\in\bM_n$, $i=1,\ldots, m$, and let $\Phi:\bM_n\to \bM_d$ be a positive linear map. Then, for  all symmetric norms and all non-decreasing geometrically convex function $g(t)$, the map
$$
(t_1,\ldots,t_m) \mapsto \left\| g\left(\Phi\left(\sum_{i=1}^m X_i^*A_i^{t_i}X_i\right)\right) \right\|
$$
is jointly log-convex on $\bR^m$.
\end{theorem}

\vskip 10pt
We will prove in the next section these two theorems. Here are some  special cases of Theorem \ref{th-norm}.

\vskip 10pt
\begin{cor}\label{cor-intro}   Let $A, Z\in\bM_n^+$.  Then, for  all symmetric norms and all non-decreasing geometrically convex function $g(t)$, 
$$
\left\| g(Z\circ I) \right\|^2
\le
\left\| g(Z\circ A) \right\| \cdot 
\left\| g(Z\circ A^{-1}) \right\| .
$$
\end{cor}

\vskip 10pt
\begin{cor}\label{cor-intro2}    Let $A_i\in\bM_n^+$ and $X_i\in\bM_n$, $i=1,\ldots, m$.  Then, for  all symmetric norms and all non-decreasing geometrically convex function $g(t)$, 
$$
\left\| g\left( \sum_{i=1}^m X_i^*X_i \right)\right\|^2
\le
\left\|  g\left(  \sum_{i=1}^m X_i^*A_iX_i \right)\right\|\cdot
\left\|  g\left(  \sum_{i=1}^m X_i^*A_i^{-1}X_i \right)\right\|.
$$
\end{cor}

\vskip 10pt
\begin{cor}\label{cor-intro3}    Let $A_i\in\bM_n^+$ and $\lambda_i>0$, $i=1,\ldots, m$, such that $\sum_{i=1}^m\lambda_i=1$. let $p> 1$ and $p^{-1}+q^{-1}=1$. Then, for  all symmetric norms and all non-decreasing geometrically convex function $g(t)$, 
$$
\left\| g\left( \sum_{i=1}^m \lambda_iA_i\right)\right\|
\le
\left\|  g\left(  I\right)\right\|^{1/q}\cdot
\left\|  g\left(  \sum_{i=1}^m \lambda_iA_i^p \right)\right\|^{1/p}.
$$
\end{cor}

 If $f(t)$ and $g(t)$ are geometrically convex then so are $f(t)+g(t)$, $\max\{f(t), g(t)\}$, $f(t)g(t)$, $e^{ f(t)}$ and $f^{\alpha}(t)$ for all $\alpha>0$. Hence the above results may be applied to a large class of functions, for instance
$$
g(t)=\sum_{k=1}^p
 c_k t^{\alpha_k}, \qquad c_k>0,\  \alpha_k \ge 0
$$
or 
$$
g(t) =\max\{ c, \beta t^{\alpha}\}, \qquad c, \alpha, \beta \ge 0.
$$
Some interesting examples of geometrically convex (also called multiplicatively convex) functions defined on a sub-interval of the positive half-line are given in \cite{Nicu}. These functions can be used to obtain exotic matrix inequalities.  A recent study \cite{BL} of a two variables log-convex functional have provided many classical and new matrix inequalities.

\section{Geodesics  and log-majorization}

The space $\bP_n$ of $n$-by-$n$ positive definite matrices is a symmetric Riemannian manifold.  There exists a unique geodesic joining  two distinct points $A,B\in\bP_n$, that can be parametrized as
\begin{equation}\label{eq1}
t\mapsto A\#_tB = A^{1/2}(A^{-1/2}BA^{-1/2})^tA^{1/2}, \qquad t\in(-\infty,\infty).
\end{equation}
In particular, the middle point between $A$ and $B$ is $A\#_{1/2}B$, the geometric mean, often merely denoted as $A\#B$. For a general $t$, especially when $t\in(0,1)$, $A\#_tB$ is a weigthed geometric mean. We refer to \cite{Bhatia} for a background on the geometric mean and  $\bP_n$.

Given $S,T\in\bM_n^+$, the weak log-majorization  relation $S\prec_{w\!\log} T$ means that
$$
\prod_{j=1}^k \lambda_j(S) \le \prod_{j=1}^k \lambda_j(T) 
$$
for all $k=1,\ldots,n$, where $\lambda_1(\cdot)\ge \ldots \ge \lambda_n(\cdot)$ stand for the eigenvalues arranged in nonincreasing order. We denote by $S^{\downarrow}$ the diagonal matrix with the eigenvalues $\lambda_1(S), \dots, \lambda_n(S)$ down to the diagonal.

\vskip 10pt
\begin{theorem}\label{thsev}   Let $A_i,B_i\in\bP_n$, $i=1,\ldots, m$ and let $\Phi:\bM_n\to \bM_d$ be a positive linear map. Then, for  all symmetric norms and all non-decreasing geometrically convex function $g(t)$, the map
$$
(t_1,\ldots,t_m) \mapsto \left\| g\left(\Phi\left(\sum_{i=1}^m A_i\#_{t_i} B_i\right)\right) \right\|
$$
is jointly log-convex on $\bR^m$.
\end{theorem}

\vskip 10pt
\begin{proof} Let $A,B\in\bP_n$ and let $\Psi:\bM_n\to\bM_d$ be a positive linear map.
We first prove the single variable case of the theorem by showing that the function
\begin{equation}\label{eqsingcase}
t\mapsto \| g(\Psi(A\#_tB))\| 
\end{equation}
is log convex on $(-\infty,\infty)$. From Ando's operator inequality
$$
\Psi(A\#B) \le \Psi(A) \# \Psi(B)
$$
and the relation $\Psi(A) \# \Psi(B)=\Psi(A)^{1/2} V\Psi(B)^{1/2}$ for some unitary $V\in\bM_d$, we infer by Horn's inequality, the weak log-majorization
$$
\Psi(A\#B) \prec_{w\!\log} \Psi(A)^{1/2\downarrow} \Psi(B)^{1/2\downarrow}
$$
Since $g(t)$ is geometrically convex, we have $g(e^{(a+b)/2})\le \sqrt{g(e^a)g(e^b)}\le (g(e^a)+g(e^b))/2$. Hence $t\mapsto g(e^t)$ is a nondecreasing convex function on $(-\infty,\infty)$. The above weak log-majorization then ensures that
$$
g(\Psi(A\#B)) \prec_{w} g(\Psi(A)^{1/2\downarrow} \Psi(B)^{1/2\downarrow})
$$
and using that $g(t)$ is geometrically convex, we infer
$$
g(\Psi(A\#B)) \prec_{w} g\left(\Psi(A)\right)^{1/2\downarrow} g\left(\Psi(B)\right)^{1/2\downarrow}.
$$
This weak majorization says that
$$
\| g(\Psi(A\#B)) \| \le \| g\left(\Psi(A)\right)^{1/2\downarrow} g\left(\Psi(B)\right)^{1/2\downarrow}\|
$$
for all symmetric norms. The Cauchy-Schwarz inequality for symmetric norms yields
$$
\| g(\Psi(A\#B)) \| \le \| g(\Psi(A))\|^{1/2} \|g(\Psi(B)\|^{1/2}.
$$
Since $A\#_{(s+t)/2} B=(A\#_sB)\#(A\#_tB)$, we get
\begin{equation}\label{eqfirst}
\| g(\Psi(A\#_{(s+t)/2}B)) \| \le \| g(\Psi(A\#_sB))\|^{1/2} \|g(\Psi(A\#_tB))\|^{1/2},
\end{equation}
for all $s,t\in(-\infty,\infty)$,
thus \eqref{eqsingcase} is a log-convex function.

We turn to the severable variables case. Let $\Phi:\bM_n\to\bM_d$ be a positive linear map, and let $A_i,B_i\in\bP_n$, $i=1,\ldots,m$. Consider the two block diagonal matices in $\bM_m(\bM_n)$,
$$
A=A_1\#_{s_1}B_1\oplus\cdots \oplus A_m\#_{s_m}B_m, \quad B=A_1\#_{t_1}B_1\oplus\cdots \oplus A_m\#_{t_m}B_m,
$$
so that
$$
A\#_{1/2} B= A_1\#_{\frac{s_1+t_1}{2}}B_1\oplus\cdots \oplus A_m\#_{\frac{s_m+t_m}{2}}B_m.
$$
Define  the positive linear map  $\Psi: \bM_m(\bM_n)\to \bM_n$,
$$
\Psi([A_{i,j}]):= \Phi\left(\sum_{i=1}^m A_{i,i}\right).
$$
From \eqref{eqfirst} with $s=0$, and $t=1$, we get
$$
 \left\| g\left(\Phi\left(\sum_{i=1}^m A_i\#_{\frac{s_i+t_i}{2}} B_i\right)\right) \right\| \le 
 \left\| g\left(\Phi\left(\sum_{i=1}^m A_i\#_{s_i} B_i\right)\right) \right\|^{1/2}
\left\| g\left(\Phi\left(\sum_{i=1}^m A_i\#_{t_i} B_i\right)\right) \right\|^{1/2}
$$
which completes the proof.
\end{proof}

\vskip 10pt
\begin{cor}\label{corsev1} Let $\varphi(t)$ be a nondecreasing function defined on $(0,\infty)$. Suppose that $\exp \varphi(t)$ is geometrically convex and let $A_i,B_i\in\bP_n$, $i=1,\ldots, m$. Then, the map
$$
(t_1,\ldots,t_m) \mapsto {\mathrm{Tr\,}}\varphi\left(\sum_{i=1}^m A_i\#_{t_i} B_i\right)
$$
is jointly convex on $\bR^m$.
\end{cor}

\vskip 10pt
\begin{proof} Let $\varphi (t)=\log g(t)$, where $g(t)$ is geometrically convex. Since $g^{\alpha}(t)$ is also geometrically convex for all $\alpha>0$,  Theorem \ref{thsev} with the normalized trace norm shows that the map
$$
(t_1,\ldots,t_m) \mapsto \frac{1}{n} {\mathrm{Tr\, }} g^{\alpha}\left(\sum_{i=1}^m A_i\#_{t_i} B_i\right)
$$
is jointly log-convex, and so is
$$
(t_1,\ldots,t_m) \mapsto \left\{\frac{1}{n} {\mathrm{Tr\, }} g^{\alpha}\left(\sum_{i=1}^m A_i\#_{t_i} B_i\right)\right\}^{1/\alpha}.
$$
Letting $\alpha\searrow 0$, we infer that the map
$$
(t_1,\ldots,t_m) \mapsto {\mathrm{det }}^{1/n} g\left(\sum_{i=1}^m A_i\#_{t_i} B_i\right)
$$
is jointly log-convex. Thus the map
$$
(t_1,\ldots,t_m) \mapsto \log {\mathrm{det\, }} g\left(\sum_{i=1}^m A_i\#_{t_i} B_i\right) = {\mathrm{Tr\,}}\varphi\left(\sum_{i=1}^m A_i\#_{t_i} B_i\right)
$$
is jointly convex.
\end{proof}

Theorem \ref{thsev} can be regarded as a generalized  H\"older inequality. This is more transparent for a single variable and pairs of commuting operators. Note that for two commuting positive definite matrices, $A\#_tB=A^{1-t}B^t$. Letting $t=q^{-1}\ (=0p^{-1}+ 1q^{-1})$ and using Theorem \ref{thsev} yields our next and last corollary.

\vskip 10pt
\begin{cor}\label{cor-holder}    Let $A_i,B_i\in\bM_n^+$ such that $A_iB_i=B_iA_i$,  $i=1,\ldots, m$. Let $p> 1$ and $p^{-1}+q^{-1}=1$. Then, for  all symmetric norms and all non-decreasing geometrically convex function $g(t)$, 
$$
\left\| g\left( \sum_{i=1}^m A_iB_i\right)\right\|
\le
\left\|  g\left(  \sum_{i=1}^m A_i^p \right)\right\|^{1/p}\cdot
\left\|  g\left(  \sum_{i=1}^m B_i^q \right)\right\|^{1/q}.
$$
\end{cor}

\vskip 10pt
Choosing $g(t)=\sinh t$, we recapture the H\"older inequality of the Abstract.

We close the paper by showing that Theorem 2.1 is equivalent to Theorem 1.2 (and similarly for  Corollary \ref{corsev1} and Theorem \ref{th-trace}). To this end, first note that by a limit argument we may assume that, in Theorem 1.2, $X_i$ and $A_i$ are invertible, $i=1,\dots,m$. Then, using the polar decomposition $X_i=U|X_i|$, observe that
$$
X_i^*A^{t_i}X_i= |X_i|(U^*AU)^{t_i}|X_i|=C\#_{t_i}D
$$
with $C=|X_i|^2$ and $D=|X_i|U^*AU|X_i|=X_i^*AX_i$.

\vskip 5pt
\noindent
Laboratoire de math\'ematiques, 
Universit\'e de Franche-Comt\'e, 
25 000 Besan\c{c}on, France.
Email: jcbourin@univ-fcomte.fr

\vskip 5pt
\noindent
College of Mathematics and Statistic Sciences, Ludong University, Yantai 264001, China.
Email: jingjing.shao86@yahoo.com

\end{document}